\documentstyle{article}
\begin{document}

\noindent
Ukrainian Mathematical Journal, Vol. 49, No. 10, 1997, 1373-1384.

\bigskip
\noindent
\centerline {\large \bf VECTOR FIELDS WITH A GIVEN SET}

\medskip
\centerline {\large \bf OF SINGULAR POINTS}

\bigskip
\centerline {\large \bf A.O.Prishlyak}

\medskip
\centerline {Kiev University, Ukraine.}

\medskip
e-mail: prish@mechmat.univ.kiev.ua

1991 MSC: 57R25, 58F25 (57R40)

\bigskip
Theorems on the existence of vector fields with given sets of Indexes of
isolated Singular points are proved for the cases of closed manifolds, pairs of
manifolds, manifolds with boundary, and gradient fields. It is proved that, on
a two-dimensional manifold, an index of an isolated Singular point of the
gradient field is not greater than one.

\bigskip
In the present paper, we consider vector fields on manifolds with isolated
singular points. In 1885, Poincare [l] proved that the sum of the indices of
the singular points of such a field on a two-dimensional manifold is equal to
the Euler characteristic of this manifold. For the $n$-dimensional case, this
fact, called the Poincare-Hopf theorem, was proved by Hopf [2] in 1926 after
partial results of Brauer and Hadamard. This theorem holds for a manifold with
boundary if the field is directed outside in any point of the boundary. It was
established that there exists a vector field without singular points on a
manifold with zero Euler characteristic. These facts are proved in [3] and [4].

The aim of the present paper is to prove the existence of a vector field with a
given set of indices satisfying the conditions of the Poincare-Hopf theorem. In
Sec. $l$, we establish the existence of such fields for closed manifolds. There
we introduce the following two operations over vector fields: introduction of a
pair of singular points and composition of singular points. These operations
are also used in the proofs of other theorems. In Sec. 2, we prove that there
exists a vector field with two sets of indices on a pair of manifolds. The case
where a vector field given on the manifold is tangent to a submanifold is
considered in particular. Manifolds with a boundary are investigated in Sec. 3.
Sections 4 and 5 deal with gradient fields for functions on manifolds. Theorem
7 is a "well-known fact," the proof of which is not yet published.

\medskip
All manifolds, functions, and vector fields considered in the present paper are
$C^{\infty}$-differentiable.

\bigskip
\noindent
{\bf 1. Singular Points of Differential Equations}

\medskip
Let $M^n$ be a smooth manifold. A vector field $v$ setting the differential
equation

\begin{equation}
\frac{dx}{dt}=v(x)
\end{equation}

\noindent
is considered as a cut of a tangent fibre bundle $TM^n$. In what follows, we
assume that the differential equation (l) is given if the vector field $v$ is
given.

\medskip
{\bf Theorem l.} Let $M^n$, $n \leq  2$, be a smooth connected manifold. Let
$\alpha _1,...,\alpha _k$, $k \leq  1$, be an integer set

$$ \sum_{i=1}^{k}{\alpha _i}=\chi (M^n) $$

\medskip
\noindent
where $\chi (M^n)$ is the Euler characteristic of the manifold $M^n$. Then
there exists a vector field $v$ on $M^n$ singular points of which are isolated
and have indices $\alpha _1,...,\alpha _k.$

\medskip
{\bf Proof.} Let us describe two operations over vector fields that allow one
to obtain a vector field with a prescribed set of singular points starting with
an arbitrary vector field.

1. Introduction of two singular points with indices +1 and -1. Let $x_0$ be a
regular point of a vector field $v_0$. According to the theorem on
rectification of a vector field [5], there exists a map $U$ at the point $x_0$,
where the vector field is constant. Thus, there exists a function $f: U
\rightarrow R^1$ without critical points in a neighborhood $U$ of the
point $x$ and such that the gradient field of this function in a proper metric
coincides with the field $v_0$. In a neighborhood $V$ of the point $x$ such
that $ \overline{V} \subset U$, we change the function $f$ by introducing a pair
ot mutually reducible critical points of adjacent indices. Then by replacing
the field $v_0$ by the gradient field of the new function in the neighborhood
$U$, we obtain a field that has two extra points of indices +1 and -1 in
comparison with the field $v_0.$

2. Replacement of two singular points of indices $\lambda _1$ and $\lambda _2$
by a singular point of index $\lambda _1+\lambda _2$. Let $x_0$ and $x_1$ be
singular points of indices $\lambda _1$ and $\lambda _2$, respectively, of a
vector field $v_0$. We select a path

$$\gamma : [0,1] \rightarrow  M^n$$

\medskip
\noindent
without self-intersections that joins points $x_0$ and $x_1$, i.e., a path
such that $\gamma (0)=x_0$ and $\gamma (l)==x_1$. Let $U$ be a neighborhood of
the path $\gamma ([0, 1])$, the closure of which contains no singular points
except $x_0$ and $x_1$. According to the Uryson lemma, there is a smooth
function $f: M^n\rightarrow R^1$ such that

\medskip
\noindent

$$f(x) = 0 \ \  for \ \ x \in  \gamma ([0,l]),$$

$$f(x) = 1   \ \ for \ \ x\in M^n\backslash U,$$

$$0 < f(x) < 1 \ \  for \ \ x\in U\backslash \gamma ([0, 1]).$$

\medskip

We consider a vector field $v_1=fv_0$. Let $g: M^n\rightarrow  M^n/\gamma ([0$,
1]) be a mapping that transforms the path $\gamma ([0,l])$ into the point $y_0$
and is a one-to-one correspondence on $M^n\backslash \gamma ([0,1])$. Since
$v_1=0$ for $x\in  \gamma ([0$, I]), the field $v_1$ induces the vector field
$v_2$ on the manifold $M^n/\gamma ([0$, 1]). In this case, the trajectories of
the vector field $v_1$ ending in (or starting from) $\gamma ([0$, 1])
correspond to trajectories ending in (or starting from) the point $y_0$. Since
the manifold $M^n/\gamma ([0,l])$ is diffeomorphic to the manifold $M^n$ the
field $v_2$ sets a field $v$ on the manifold $M^n$. By virtue of the
Poincare-Hopf theorem, the sum of the indices of $v$, just as $v_0$, is the
Euler characteristic $\chi (M^n)$. Thus, $y_0$ is a singular point of the index
$\lambda _1+\lambda _2$ of the vector field $v.$

Let $f: M^n \rightarrow R^1$ be a Morse function on the manifold $M^n$. We
transform this function by introducing pairs of mutually reducible critical
points in such a way that the number of critical points is no less than . Then
the gradient field of the function $f$ consists of singular points of indices
$+1$ and -1 corresponding to critical points of the function $f$. We decompose
the set of these critical points into groups such that the sum of the indices
of the $i-th$ group is $\alpha _i$ and then apply operation 2 to the points of
every group. As a result, we get the desired field $v.$

\bigskip
\noindent
{\bf 2. Differential Equations on Pairs of Manifolds}

\medskip
{\bf Definition 1.} Assume that $M^n$ is a smooth manifold, $N^k$ is its
submanifold, $\rho $ is a Riemann metric, and $v$ is a vector field on the
manifold $M^n.We$ say that a vector field $u$ on the submanifold $N^k$ is
induced by the
vector field $v$ in the Riemann metric $\rho $ if, for any point $x\in N^k$ ,
the vector $u(x)\in T_xN^k$ is the projection orthogonal in the metric $\rho $
of the vector $v(x)\in T_xM^n$ to the subspace $T_xN^k.$

If $v$ is the gradient field of the function $f$ on the manifold $M^n$ then the
induced vector field $u$  on the submanifold $N^k$ coincides with grad $g$,
where $g$ is the restriction of the function $f$ to the submanifold $N^k.$

\medskip
{\bf Theorem 2.} Assume that $N^k$ is a submanifold of a smooth manifold $M^n$,
$\rho $ is a metric on the manifold $M^n$ and $\alpha _1,...$, $\alpha _p
(p\leq 1)$, $\hbox{$\beta $}_1,...$, $\hbox{$\beta $}_s (s \leq  1)$ are
integer sets such that

$$ \sum_{i=1}^{p}{\alpha _i}=\chi (M^n), \ \
 \sum_{i=1}^{s}{\beta  _i}=\chi (N^k), $$

\medskip
\noindent
where $\chi (M^n)$ and $\chi (N^k)$ are the Euler characteristics of the
manifolds $M^n$ and $N^k$ respectively. If $n-k \leq 2$, then on the manifold
$M^n$ there exists a vector field $v$ with singular points of indices
$\alpha _1,...$, $\alpha _p$ such that the vector field $u$ on the submanifold
$N$ induced by the vector field $v$ in the metric $\rho $ has singular points
of indices $\hbox{$\beta $}_1,...$, $\hbox{$\beta $}_s.$

\medskip
{\bf Proof.} By using Theorem 1, we construct a vector field $v$ on $M^n$ with
singular points of indices $\alpha _1,...$, $\alpha _p$. Then we consider this
vector field as a section of the tangent fiber bundle and present it in a
general position with respect to the submanifold $N^k$ located in the zero
section of the tangent fiber bundle $TM^n$. The obtained vector field $v_1$ has
singular points with the same indices $\alpha _1,...,\alpha _p$; these points
are not located on the submanifold $N^k$. Let A be a tubular neighborhood of
$N^k$ that contains no singular points of the vector field $v_1$. We consider
the restriction of the tangent fiber bundle $TM^n$ to the manifold $N^k$. The
obtained fiber bundle $\xi $ has dimension $n+k$ as a manifold. The tangent
fiber bundle $TN^k$ is a subbundle of $\xi $ and has dimension 2k as a
manifold. The restriction of the vector field $v$ to the submanifold $N^k$
considered as a section of the fiber bundle $\xi $ is of dimension $k$. In
general position, its intersection with the tangent fiber bundle $TN^k$ is a
submanifold $L$ of dimension $2k -n.$

Let $u_1$ be a vector field on the manifold $N^k$ induced by the vector field
$v_1$. Evidently $u_1(x)= v_1(x)$ if and only if the point $x$ belongs to the
submanifold $L$, and $u_1(x)=0$ if and only if the vector $v_1(x)$ is
perpendicular to the submanifold $N^k$ with respect to the metric $\rho $.
Thus, singular points of the vector field $u_1$ are not on the submanifold L.
By using the vector field $u_1$, a vector field $u$ on the manifold $N^k$ can
be constructed such that $u(x) == u_1(x)$ for all points $x$ of some
neighborhood of the submanifold L. In fact, two singular points of indices $+1$
and -1 can be introduced in a neighborhood $V$ of any nonsingular point $x_0$
if $V$ contains no singular points and the closure of $V$ does not intersect
the submanifold L. Since the codimension of the submanifold $L$ in the manifold
$N$ is at least two, we can join any two singular points in a way that does not
intersect $L$ and choose a neighborhood $W$ of this path such that the closure
of the neighborhood does not intersect the submanifold $L$ and contains no
other singular points. Let us replace, as was done in Theorem 1, these two
singular points with indices $\lambda _1$ and $\lambda _2$ by a singular point
of index $\lambda _1+ \lambda _2$. In this case, the vector field is changed
only on the set $W.$

Let us consider the vector field $u-u_1$ on the manifold N. We extend this
vector field to a vector field $w$ in a tubular neighborhood A of the
submanifold $N^k$ as follows: the coordinates of the vector $w(x)$ are equal to
the coordinates of the vector $u(x_0)-u_1(x_0)$ in some map if the point $x$ is
located in a fiber bundle of the tubular neighborhood over the point $x_0$. Let
us take the tubular neighborhood A small enough so that vectors $w(x)$ will not
be parallel to the vectors $v(x)$ if $w(x) \neq  0.$

By the Uryson lemma, there is a smooth function $f: M^n \rightarrow R^1$
such that

\medskip

$$f(x) =1 \ \ for \ \ x \in  N^k$$

$$f(x) = 0 \ \   for \ \ x \in  M^n\backslash A,$$

$$0 <f(x) < 1 \ \   for \ \ x \in  A\backslash N^k.$$

\medskip
\noindent
We consider a vector field $v = v_1+f w$. According to the construction, it has
the same set of singular points on $M^n$ as the field $v_1$ does. The field
induced by $v_1$ on the manifold $N^k$ coincides with the field $u = u_1+1 (u -
u_1)= v_1 +f w$. Theorem 2 is proved.

Now let us consider the case of vector fields on a pair of manifolds $(M^n$,
$N^k)$ such that the manifold $N^k$ is imbedded into the manifold $M^n$, where
the restriction of a field $v$, given on the manifold $M^n$, to the manifold
$N^k$ is a vector field tangent to $N$ , $i.e.$, the vectors induced by the
field $v$ coincide with the corresponding vectors of the field $v$. It is
evident that the induced vector field in this case does not depend on the
metric $\rho $ on the manifold $M^n$ and every singular point of the vector
field $v$, located on the submanifold $N^k$, is a singular point of the vector
field $v$ set on the manifold $M^n.$

\medskip
{\bf Theorem 3.} Let $N^k$ be a submanifold of a smooth manifold $M^n$, $n
-k\leq 1$,  and let $\alpha _1,...,\alpha _p$ and $\hbox{$\beta $}_1,...$,
$\hbox{$\beta $}_s$ be integer sets such that

$$ \sum_{i=1}^{p}{\alpha _i}=\chi (M^n), \ \
 \sum_{i=1}^{s}{\beta  _i}=\chi (N^k), 1\leq s \leq  p.$$

\medskip
\noindent
Then, there exists a vector field $v$ on the manifold $M^n$ with singular
points of indices $\alpha _1,...,\alpha _p$ tangent to the submanifold $N^k$,
and such that the vector field $v$ on the submanifold $N^k$ has singular points
of indices $\hbox{$\beta $}_1,...$, $\hbox{$\beta $}_s.$

\medskip
{\bf Proof.} Let us construct, by analogy with the proof of Theorem 1, a vector field
$u_1$ on the submanifold $N^k$ with singular points of indices
$\hbox{$\beta $}_1,...$, $\hbox{$\beta $}_s$. We select on the manifold $M^n$ a
vector field $u_2$ which is normal to the submanifold $N^k$ and nonzero in the
singular points of the vector field $u_1$ . Let $U$ be a tubular neighborhood
of the submanifold $N$ . According to theUryson lemma, there exists a smooth
function $f$ on the manifold $M^n$ such that

\medskip

$$f(x) =0 \ \  for \ \ x \in  N^k,$$

$$f(x) =1 \ \  for \ \ x \in  M^n\backslash U,$$

$$0 <f(x) < 1 \ \   for \ \ x \in  U\backslash N^k.$$

\medskip
\noindent
Let  $u = (l -f) u_1+f\ u_2$ be a vector field on $U$, where $u_1$ and $u_2$
are vector fields on $U$, the vectors of which in every fiber bundle on a
tubular neighborhood have the same coordinates in some map as the corresponding
vectors on $N^k$ do. Let us arbitrarily extend the field $u$ onto the manifold
$M^n$. When considering $u$ as a section of the tangent fiber bundle,  we  lead
it  to  the  general  position  with  the  zero
section and preserve it without changes on the set U. We introduce pairs of
singular points with indices +1 and -1 and add singular points along paths
the inner parts of which do not intersect the submanifold $N^k$. The addition
of two points located on the submanifold is not assumed. The vector field thus
constructed is a tangent field for the submanifold $N^k$ and satisfies all the
conditions of the theorem.

\medskip
{\bf Remark 1.} For every singular point on the submanifold $N^k$, any two
integer numbers can be the indices of the vector field $v$ and the indices of
the vector field $u.$

\medskip
{\bf Remark 2.} If $n- k>2$ or $n - k = 1$ and the submanifold $N^k$ does not
divide the manifold $M^n$, then the theorem also holds for $s=0$. That is, if
$\chi (N^k)==0$, then there exists a vector field $v$ on the manifold $M^n$
that is tangent to the submanifold $N^k$ has no singular points on $N^k$, and
possesses the given set of singular points $\alpha _1,...$, $\alpha _p$ on the
submanifold $M^n (\Sigma ^p_{i=1} \alpha _i= \chi (M^n)).$

If the submanifold $N^k$ decomposes the manifold $M^n$ into two submanifolds
$M_1$ and $M_2$ with the boundary $\partial M_1 = \partial M_2 = N^k$, then the
problem of existence of a vector field $v$ on the manifold $M^n$ with a given
set of indices and tangent to the submanifold $N^k$ and without singular points
on $N$ is equivalent to the problem of existence on a manifold with a boundary
of a vector field tangent to the boundary and the singular points of which are
not located on the boundary and have the preassigned set of indices.

\bigskip
\noindent
{\bf 3. Differential Equations on Manifolds with Boundaries}
\medskip

{\bf Proposition 1.} Let $M^n$ be a smooth manifold with the boundary $N$,
$\chi (N)=0$. On the manifold $M$, there exists a vector field tangent to the
manifold $N$, and all its singular points are inner and have indices
$\alpha _1,...$, $\alpha _p$ if and only if

$$ \sum_{i=1}^{k}{\alpha _i}=\chi (M^n), $$

\noindent
where $\chi (M^n)$ is the Euler characteristic of the manifold $M$ with the
boundary $N.$

\medskip
{\bf Proof.} Since the Euler characteristic of the manifold $N$ is 0, there exists a
vector field $u$ on $N$ without singular points. Let $w$ be a vector field on
$N$, normal to $N$, without singular points, and directed outside of the
manifold $M$, We extend this field onto the collar $N\times [0$, 1 ] according
to the formula

$$v(x, t) = tu(x) + (1 - t)w(x),$$

\medskip
\noindent
where the point $x\in N$, $t\in  [0,1]$. We extend the field $v$ onto the
manifold $M^n$ in such a manner that all singular points are isolated. Then the
sum of the indices of these singular points is equal to the Euler
characteristic $\chi (M^n)$. By using the processes of the introduction of new
pairs of singular points of indices $+ 1$ and -1 and the composition of
singular points, we obtain the desired vector field.

Let us consider an arbitrary vector field $v$ on the manifold $M^n$ with the
boundary $\partial M^n = N$. We assume that the vector field $v$ has no
singular point on the boundary $N$, Let $u$ be a vector field induced by the
vector field $v$ on the boundary N. We consider a manifold $M'$ obtained from
the manifold $M^n$ by gluing the collar $N\times [0$, 1]; to do this, we
identify points $x \in  N = \partial M^n$ with points $(x$, $0) \in
N\times [0$, 1]. We construct a vector field $w$ on $M'$ which is an extension
of the field $v$. Let a be a vector field given on $N\times \{l\}$, orthogonal
to $M'$ directed outside of the manifold $M'$ and such that

$$|a| > |v(x,0)|$$

\medskip
\noindent
for any $x \in  N$. We extend the field $v$ onto the collar $N \times  [0$, 1]
according to the relation

$$v(x,t) = ta+ (1-t)v(x,0).$$

\medskip
\noindent
By smoothing this field, we get a vector field $w$. It is evident that, for a
point $x \in  N$, there exists $t\in (0$, 1) such that $w(x$, $t) = 0$ if and
only if $u(x) = 0$ and the vector $v(x$, 0) is directed inside of the manifold
$M^n$. The indices of the corresponding singular points have the same absolute
value and different signs. We denote by $\delta _+(v)$ and $\delta _- (v)$,
respectively, the sum of the indices of singular points of the vector field $u$
where the vector field $v$ is directed inside of and out of the manifold $M^n$.
Then the following equalities hold:

$$\chi (M^n) = \sum_{}^{}{ind\ v} - \delta _+(v),$$

$$\chi (M^n) =  \sum_{}^{}{ind\ v} + \delta _-(v)\ \  for \ even \ n,$$

$$\chi (M^n) = - \sum_{}^{}{ind \ v} - \delta _-(v)\ \  for \ odd \ n.$$

\medskip
Here $\chi (M^n)=\chi (M')= \sum_{}^{}{ind\ w}$  is the Euler characteristic of
the manifold $M^n$ with the boundary $N.$

\medskip
{\bf Theorem 4.} Let $M^n$ be a smooth manifold with the boundary $N$, let
$(\alpha _1,...,\alpha _s$, $\hbox{$\beta $}_1,...,\hbox{$\beta $}_k$,
$\hbox{$\beta $}_{k+1},...,\hbox{$\beta $}_{p+1})$ be integer sets. There
exists a vector field $v$ with inner singular points of indices
$\alpha _1,...$, $\alpha _s$ and with the induced vector field $u$ on the
manifold $N$ such that the vector field $v$ is directed inside of the manifold
$M^n$ at the singular points of the field $u$ with indices
$\hbox{$\beta $}_1,..$. , $\hbox{$\beta $}_k$ and outside of this manifold at
the singular points of the field $u$ with indices $\hbox{$\beta $}_{k+1},...$,
$\hbox{$\beta $}_p$ if and only if

\medskip
$$ \sum_{i=1}^{p}{\beta _i}=\chi (N)$$

$$\chi (M^n) =\sum_{i=1}^{s}{\alpha _i} - \sum_{i=1}^{k}{\beta _i} \ for\ even\ n,$$

$$\chi (M^n) =\sum_{i=1}^{s}{\alpha _i} - \sum_{i=k+1}^{p}{\beta _i} \ for\ odd\ n.$$

\medskip
{\bf Proof.} The previous reasoning demonstrates that the conditions of the
theorem are satisfied for any vector field $v$. Let us prove the inverse
statement $i.e.$, if the sets of indices satisfy the conditions of the theorem,
then there exists a vector field with given sets of indices. In fact, let $u_0$
be an arbitrary vector field tangent to the manifold $N=\partial M^n$, singular
points of which have indices $\hbox{$\beta $}_1,...,\hbox{$\beta $}_k$,
$\hbox{$\beta $}_{k+1},...$, $\hbox{$\beta $}_p$. Such a vector field exists by
virtue of Theorem 1. Let $u_1$ be a vector field given on the boundary $N$,
orthogonal to $N$, directed inside of the manifold $M^n$ in singular points
with indices $\hbox{$\beta $}_1,...$, $\hbox{$\beta $}_k$ and outside of $M^n$
in singular points with indices $\hbox{$\beta $}_{k+1},...$,
$\hbox{$\beta $}_p$. We consider a field

$$u_2 = u_0+u_1.$$

\noindent
Let us arbitrarily extend this field up to the field $v_0$ on the whole of the
manifold $M^n$. Then the field $v_0$ induces a vector field $u_2$ on the
manifold $N$ with the desired set of indices. The conditions of the theorem are
satisfied for this field. By applying the processes of introduction of a pair
of singular points of indices $+1$ and -1 and composition of singular points
described in Theorem 1 to the field $v_0$, we get the desired vector field $v.$

\bigskip
\noindent
{\bf 4. Gradient Field of a Smooth Function}
\medskip

{\bf Definition 2}. Let $v$ be a vector field on a manifold $M^n$ that has only
isolated singular points. Let us consider an oriented graph, the vertices
$\alpha _i$ of which correspond one to one with singular points $x_i$ of the
vector field $v$, and two vertices are joined by an arc if there exists an
integral trajectory that starts and ends at the corresponding points of the
vector field. We call such a graph the graph $G(v)$ of the vector field $v.$

\medskip
{\bf Definition 3}. The contour of a graph $G$ is a sequence $S = ((\alpha _0$,
$\gamma _1$, $\alpha _1$, $\gamma _2, ..., \alpha _{n-1}$,
$\gamma _n$, $\alpha _n)$ of its vertices $\alpha _i$ and arcs $\gamma _i$
alternating in such a manner that $\alpha _i-_1$ is the beginning and
$\alpha _i$ is the end of the arc $\gamma _i$, $i=\overline{1,n}$, and
$\alpha _0=\alpha _n.$

\medskip
{\bf Theorem 5.} Let $v$ be a smooth vector field on a smooth manifold $M^n$
that satisfies the following conditions:

\medskip
(i) all singular points are isolated; limit set of each trajectories consist of
singular points; for every singular point $x_i$ there is a neighborhood $U_i$
and a smooth function $f_i$ determined on $U_i$ and such that $v$ is a gradient
field of the function $f_i$ in some metric $\rho _i$ on $U_i;$

\medskip
(ii) the graph $g(v)$ of the vector field $v$ has no contours.

\medskip
Then there is a Riemann metric $\rho $ on the manifold $M^n$ and a function $f:
M^n\rightarrow R^1$  such that grad $(f) = v.$

\medskip
{\bf Proof.} Let $\rho _0$ be a fixed metric on the manifold $M^n$. Since the
graph of the vector field has no contours, the singular points $x_i$ can be
ordered in such a manner that if $i<j$, then the graph $G(v)$ does not have a
path originating at the vertex $\alpha _i$ and ending at the vertex
$\alpha _j$. Without any loss of generality, we assume that any neighbornood
$U_i$ is homeomorphic to an open disk  and $U_i\bigcap U_j= \oslash $ , $i\neq j$. We
define the function $f$ on $U_i$ in such a way that $f(x_i)=i$, $f(y)=f_i(y)
+i-f_i(x_i)$, $y\in U_i$. We here take $U_i$ sufficiently small such that

$$|f(y) - f(x_i)| < \frac{1}{3},\ \  y\in U_i.$$

We define $f$ on $M^n\backslash \bigcup_i U_i$. Let $x\in M^n\backslash \bigcup_i
U_i$, $\gamma _x$ be a trajectory of the vector field $v$ passing through the
point $x$. Let $x_i$ and $x_j$ be points where the trajectory $\gamma _x$
starts and ends, $y_i=\gamma _x$\`E$\partial U_i$,
$y_j=\gamma _x$\`E$\partial U_j$ and $y_i=\gamma _x(t_i)$,
$y_j=\gamma _x(t_j)$. We set

$$f(x)=f(y_i)+ \frac{S_{\gamma} (t_i,0)}{S_{\gamma} (t_i,t_j)}(f(y_i) - f(y_j)),$$

\medskip
\noindent
where $S_\gamma (t_i,t)$ is the length of the trajectory $\gamma _x$ in the
metric $\rho _0$ between $t_i$ and $t$. Thus, the function $f$ in-creases from
$f(y_i)$ to $f(y_j)$ along $\gamma _x$ in proportion to the length of the arc
$\gamma _x$. Let us smooth the function $f$  on
the boundary $\partial U_i$ as was done in [4]. The metric $\rho $ can be taken
as follows: for any point $x$, $x^1x_i$ we choose a coordinate system $x^1$,
$x^2,...$, $x^n$, we directs $x^1$ along the integral trajectory passing
through the point $x$, and we take $x^2$ ,..., $x^n$ on the level surface of
the function $f$ that passes through the point $x$. The scalar product at a
point $x$ is taken to be in proportion to the standard one, namely,

\medskip
$$\rho (x,y)=\frac{|v(x)|}{|\partial f(x)/\partial x|}\sum_{i=1}^{n}{x^iy^i}.$$

\medskip
\noindent
By using a partition of the unit, we glue this metric with the metrics
$\rho _i$ given on $U_i$. The obtained metric is the desired one.

\medskip
{\bf Definition 4.} A smooth function $f: M^n \rightarrow R^1$  is called
minimal if any other smooth function $g: M^n \rightarrow R^1$   has no fewer
critical points than $f$ does. We denote it by $q(M^n).$

\medskip
{\bf Definition 5.} An integer set $\alpha _1,...,\alpha _k$ is called
admissible for a manifold $M^n$ if $\alpha _1 = 1$, $\alpha _k = (-1)^n$, and
$\sum_{i=1}^{n}{\alpha _i}=\chi (M^n).$

\medskip
{\bf Definition 6.} An integer set $\alpha _1,...,\alpha _k$ is called
realizable by a smooth function $f$ if the gradient field of the function $f$
in some metric $\rho $ has $k$ isolated singular points and their indices are
$\alpha _1,...,\alpha _k.$

Evidently, if $\alpha _1,...,\alpha _k$ are indices of singular points
$x_1,..$. , $x_k$ of a gradient vector field of a function $f$ and
$f(x_i)\leq f(x_j)$ if $i\leq j$, then the set $\alpha _1,...,\alpha _k$
is admissible.

\medskip
{\bf Theorem 6}. Let $M^n$ be a smooth manifold, $n\leq 4$. Let
$\alpha _1,...,\alpha _k$ be an admissible set. Then there exists a smooth
function $f$ realizing the set if and only if $k> q(M^n).$

\medskip
{\bf Proof.} Let $f$ be a minimal function. We construct a function $f_1$ having the
same number of critical points as $f$ and such that its gradient field has an
integral trajectory starting at the point of minimum and ending at the point of
maximum of $f_1$. Let $y_1$, $y_2, ..., y_s$ be critical values of the
function $f$ such that $y_i<y_j$ if $i<j$, and $x_0$ and $x_s$ are points of
minimum and maximum of the function $f$. We set

$$y_{i+1/2}=\frac{y_i+y_{i+1}}{2}.$$

\medskip
Let $z_1\in  f^{-1}(y_1)$ be a regular point of the mapping $f$, and let
$\gamma _1$ be the trajectory of the gradient field that passes through the
point $z_1$. It is evident that $\gamma _1$ starts at the point of the minimum
of the function $f$. If $\gamma _1$ ends at the point of the maximum, then
$f_1=f$. Assume that the trajectory $\gamma _1$ is ending at a critical point
$x_i$, $f(x_i) =y_i$, $i<s$. Let $z_i\in  f^{-1}(y_i),$

$$p_i=\gamma _1 \bigcap f^{-1}(y_{i-1/2}),$$

$$q_i=\gamma _i \bigcap f^{-1}(y_{i-1/2}).$$

\medskip
\noindent
Here, $\gamma _i$ is the integral trajectory passing through the point $z_i$.
The point $z_i$ is chosen in such a way that the points $p_i$ and $q_i$ are in
the same connected component of the submanifold $f^{-1}(y_{i-1/2})$. Then
there exists an ambient isotopy of the manifold $f^{-1}(y_{i-1/2})$ that
transfers the point $q_i$ into the point $p_i$. This isotopy sets a vector
field (that is a set of integral trajectories) on the submanifold
$f^{-1}([y_{i-1/2}-\hbox{$\varepsilon $}$,
$y_{i-1/2}+\hbox{$\varepsilon $}])$ for $\hbox{$\varepsilon $}>0$ sufficiently
small. By smoothing the vector field on the boundary
$f^{-1}(y_{i-1/2}-\hbox{$\varepsilon $})$ and
$f^{-1}(y_{i-1/2}+\hbox{$\varepsilon $})$ we obtain that the integral
trajectory $\gamma _1$ passes through the point $z_i$. Applying the same
reasoning to all points $z_j$, $f(z_j)<y_s$, we construct a vector field, the
trajectory $\gamma _1$ of which starts at the point $z_0$ and ends at the point
$z_s$. By applying Theorem 5 to this vector field, we construct a function $f$,
the gradient field of which has trajectories starting and ending at $z_0$ and
$z_s$ respectively, and its set of critical points coincides with that of the
function $f.$

\medskip
Evidently, the integral trajectories sufficiently close to the trajectory
$\gamma _1$ also start and end at the points $z_0$ and $z_s$ respectively. By
introducing pairs of singular points of indices $+1$ and -1 along these
trajectories and composing singular points on the corresponding critical levels
as was done in [6], we get the desired function $f.$

\bigskip
\noindent
{\bf 5. Gradient Fields on Two-Dimensional Manifolds}

\medskip
{\bf Theorem 7.} Let $y_0$ be a singular point of a gradient vector field $v$
of the function $f$ on the manifold $M^2$. Then its index is less than two.

\medskip
{\bf Proof.} We take a neighborhood $U$ of the point $y_0$ and coordinates
$(u,v)$ in it in such a manner that $y_0$ is the origin and there are no other
critical points in $B = \{(u,v): u^2 + v^2 \leq  1\}$ except the point $y_0$,
On the boundary of the circle $B^2$, we introduce a parametrization $\partial
B^2 = S^1=\{(u,v): u^2+v^2=1\}=\{(u,v): u = \cos \ t$, $v = \sin \ t,\ 0 \leq  t
\leq  2\pi \}$. We denote by $x_i$ and $\tau (t)$, respectively, a point on
the circumference $S^1$ with coordinates $(\cos \ t$, $\sin \ t)$ and the
vector tangent to $S^1$ at the point $x_t$

$$\tau (t)=\{-\sin \ t, \cos \ t\}.$$

\medskip
\noindent
We also denote by $\alpha (t)$ a continuous function of an angle put
counterclockwise between the vector $\tau (t)$ and the vector of the field $v$
at a point $x_t$. By virtue of the index of a singular point on a
two-dimensional manifold, $\alpha (2\pi )=\alpha (0)+2(k-1)\pi $, where $k$ is
the index of the singular point $y_0.$

\medskip
{\bf Definition 7}. We say that a function $\alpha (t)$ at a point $t_0$,
increasing (or decreasing), passes a level $\alpha _0$ if
$\alpha (t_0)=\alpha _0$ and there is a neighborhood $V = (t_0-
\hbox{$\varepsilon $}_1$, $t_0 +\hbox{$\varepsilon $}_2)$ of the point $t_0$
(here, $\hbox{$\varepsilon $}_1>0$ and $\hbox{$\varepsilon $}_2 > 0)$ sach that
the function $\alpha $ is monotonically nondecreasing (or nonincreasing) in the
interval $V$ and

$$ \alpha (t_0- \varepsilon _1) < \alpha _0 < \alpha (t_0
+\varepsilon _2) \ \ (or \ \alpha (t_0- \varepsilon _1) >
\alpha _0 > \alpha (t_0 +\varepsilon _2)).$$

\medskip
{\bf Lemma 1}. Assume that a trajectory $\gamma (s)$ of a vector field $v$
passes through a point $x_{t_0}$ of a circumference $S^1$ such that at the point $t_0$
the function $\alpha (t)$, increasing, passes the level $\pi n (n \in Z)$ .
Then this trajectory is situated in the circle $B^2$ locally in a neighborhood
of the point $x_{t_0}$.

\medskip
{\bf Proof.} Assume the contrary. Let the trajectory $\gamma (s)$ pass through
the point $x_{t_0}$ (i.e., $\gamma (s_0) = x_{t_0}$) and leave the circle
$B^2_{1+\varepsilon }= \{(u,v): u^2+v^2 \leq l+\varepsilon \}$ for some $\varepsilon > 0$. For
definiteness, we assume that $n=0$ and $\gamma (s)$ leaves the circle  as the
parameter $s$ increases. We take $\varepsilon $ sufficiently small
(this can be done because of the continuity of the vector field $v$) for the
inequality

$$ \alpha _{\varepsilon }(t_0- \varepsilon _1) < 0 <
\alpha _{\varepsilon }(t_0 +\varepsilon _2) $$

\medskip
\noindent
to be satisfied, where $\alpha _{\varepsilon }$ is the function of the
angle between the tangent vector to the boundary $\partial B^2_{1+\varepsilon }
= S^1_{1+\varepsilon } $ and the corresponding vector of the field $v$, the
circle $B^2_{1+\varepsilon }$ contains no other singular
points except $y_0$, and the trajectory passing through a point
$x_{t_0+\varepsilon _2}$ intersects $S^1_{1+\varepsilon }$ at
the point $y$ with the parameter $t<t_0+{\varepsilon }_2$ (for this
purpose, ${\varepsilon }_2$ should satisfy the inequality
$\alpha (t_0+{\varepsilon }_2) < \pi /2)$. Then there exists
$t_1\in (t_0-{\varepsilon }_1,\ t_0+{\varepsilon }_2)$ such that
$\alpha _{\varepsilon }(t_1)=0$. We denote this point on the
circumference  $S^1_{1+\varepsilon }$ by $y_1$. Similarly, for any ${\varepsilon }_i$,
$0<{\varepsilon }_i<{\varepsilon }$, there is $t_i$ such that
$\alpha_{\varepsilon _i}(t_i)=0$
and $t_i$ continuously depends on $\hbox{$\varepsilon $}_i$. We consider a path
$\hbox{$\beta $}$ that consists of the following three parts:

\medskip
(i) the point $x_{t_i}$ on  $S^1_{\varepsilon _i}$($\alpha_{\varepsilon _i}(t_i)=0$);

\medskip
(ii) the arc of the circumference $S^1_{\varepsilon }$ from the point $y_1$ to the point $y$;

\medskip
(iii) the arc of the trajectory of the vector field $v$ from the point $y$ to
the point $x_{t_0+\varepsilon _2}$.

\medskip
Then the trajectory $y$ leaves the circle $B^2_{1+\varepsilon }$ and must pass
through a point of the path $\hbox{$\beta $}$. This is impossible since all
vectors of the field $v$ are directed into $B^2_{1+\varepsilon }$ at points of
the path $\hbox{$\beta $}$. The contradiction obtained proves Lemma 1.

One can similarly prove the following fact. If a function $\alpha (t)$,
decreasing, passes through a level $\pi n$, $n\in Z$, at a point $t_0$, then,
locally in a neighborhood of the point $x_{t_0}$, the integral trajectory passing
through this point does not intersect the interior of the circle $B^2.$

Let us prove the following statement: if $y_0$ is a singular point of index two
and there are only two points $x_1$ and $x_2$ on the circumference $S^1$ where
the function $\alpha $, increasing, passes through the level $\pi n$, $n
\in Z$, then trajectories of the vector field $v$ passing through the points
$x_1$ and $x_2$ both start and end at the singular point $y_0$. In fact, let
$\gamma _1(s)$ and $\gamma _2(s)$ be respectively trajectories passing through
$x_1$ and $x_2$. Let $\hbox{$\beta $}_1$ be an arc of the circumference $S^1$
between the points $x_1$ and $x_2$, in the points of which the vectors of the
field $v$ are directed inside $B^2$ or are tangent to the circumference $S^1$,
and let $\hbox{$\beta $}_2$ be an arc of the circumference $S^1$ between the
points $x_2$ and $x_1$, in the points of which the vectors of the field $v$ are
directed outside the circle $B^2$ or are tangent to the circumference $S^1.$

If the trajectory $\gamma _1$ leaves the circle $B^2$ as the parameter $s$
increases (or, similarly, the parameter $s$ decreases), then it intersects the
arc $\hbox{$\beta $}_2$ in a point $x_3$. Thus, this trajectory decomposes the
circle $B^2$ into two parts $A_1$ and $A_2 y$ the arc between points $x_1$ and
$x_3.$

Assume that the points $y_0$ and $x_2$ are located in the same part $A_1$.
Since the trajectory $\gamma _1$ cannot enter the interior of the domain $A_2$
through an arc of the circumference $S^1$ between points $x_3$ and $x_1$, it
has to start in a singular point inside of the domain $A_2$. This is impossible
because the circle $B$ contains no other singular points except $y_0.$

Let the points $y_0$ and $x_2$ be located in distinct parts $A_1 i$ and $A_2$
and let the trajectory $\gamma _1$ start at the point $y_0$. We similarly
consider a trajectory $\gamma _2$. One can prove that it starts or ends at the
singular point $y_0$. However, if the points $y_0$ and $x_2$ are located in
distinct domains $A_1$ and $A_2$, then the trajectories $\gamma _1$ and
$\gamma _2$ have to intersect one another, which is impossible.

Thus, both trajectories $\gamma _1$ and $\gamma _2$ are completely located in
the circle $B^2$ and start and end at the singular point $y_0$. However, if
$\gamma _1$ and $\gamma _2$ are trajectories of the gradient field of the
function $f$, then this function strictly increases along these trajectories,
which is impossible because

\medskip

$$f(y_0) = f(\gamma _1(-\infty )) = f(\gamma _1(+\infty )) =
f(\gamma _2(-\infty )) = f(\gamma _2(+\infty ))$$

\medskip
Now let us consider the case where there are more than two points where the
function $\alpha (t)$ passes through levels $\pi n$, $n\in Z$. Note that the
difference between the number of points where the function $\alpha (t)$
increases and the number of points where it decreases when passing the levels
$\pi n, n\in Z$, is $2(k-l)$, where $k$ is the index of the singular point
$y_0.$

Let $x_i$ be a point on the circumference $S^1$, where the function
$\alpha (t)$, increasing, passes a level $\pi n$, $n\in $\`U, and let
$\gamma _i(s)$ be the integral trajectory passing through the point $x_i$. An
integral trajectory of a gradient field cannot start and end at the same
singular point $y_0$. Hence, there is a point $y_i$ where the trajectory
$\gamma _i(s)$ leaves the circle $B^2$. Then the arc of the trajectory
$\gamma _i(s)$ between the points $x_i$ and $y_i$ decomposes the circle $B^2$
into two parts. We denote the part containing the point $y_0$ by A. The
boundary of the domain A at the point $x_i$ must be smooth. If this is not the
case, we take another direction along the trajectory $\gamma _i(s)$ from the
point $x_i$. Let us smooth out the boundary $\partial A$ at the point $y_i$ and
denote the obtained curve by $S_1^1$, and the corresponding domain bounded by this
curve by $B_1^2$.

The pair $(B_1^2,S_1^1)$ is evidently diffeomorphic to the pair $(B^2$, $S^1)$. It is
possible to introduce coordinates $(u^2, v^1)$ on it such that $S_1^1$ is the unit
circle relative to these coordinates. When considering a parametrization of
this circle and introducing the function of the angle $\alpha ^1$ between a
vector tangent to $S^1_1$ and the vector field $v$, one can see that this function has at
least two points less where it passes through the levels $\pi n$, in comparison
with the function $\alpha $ (at the point $x_i$ the function $\alpha ^1$ does
not pass the levels $\pi n$; there are no points of passage through the levels
$\pi n$ on the arc of the circle $S_1^1$ between the points $x_i$ and $y_i).$

In turn, the pair $(B^2_1, S_1^1)$ is diffeomorphic to the pair $(B^2_2, S^1_2)$, for which the function
$\alpha ^2$ has at least two points less where it passes through the levels
$\pi n$, in comparison with the function $\alpha ^1$. We repeat this process
until there are two points left where a function $\alpha ^i$ passes through the
levels $\pi n$. Since the trajectories passing through these points start and
end at the singular point $y_0$, we have proved the following fact: There is no
smooth function $f$ on a two-dimensional closed manifold $M^n$ the gradient
field of which has a singular point of index $k = 2.$

For the case $k > 2$, one can demonstrate that there exists a trajectory of the
vector field originating and ending at the singular point $y_0$, which is
impossible for the gradient field of a function $f$. Theorem 7 is completely
proved.

\bigskip

\begin{center}
{\bf REFERENCES}
\end{center}
\medskip

\noindent
1. H. Poincare, On Curves Determined by Differential Equations [Russian
translation], OGIZ, Moscow-Leningrad (1947).

\noindent
2. H. Hopf, ``Vektorfelder in $n$-dimensionalen Mannigfaltigkeiten," Math.
Ann,, {\bf 96,} 209-221 (1926).

\noindent
3. J. W. Milnor, Topology from the Differentiate Viewpoint, The University
Press of Virginia, Charlottesville (1965).

\noindent
4. M. W. Hirsch, Differential Topology, Springer-Verlag, New
York-Heidelberg-Berlin (1976).

\noindent
5. V. I. Arnol'd, Ordinary Differential Equations [in Russian], Nauka, Moscow
(1971).

\noindent
6. S. Smale, "On gradient dynamical systems," Ann. Math., 74, No. 1, 199-206
(1961).

\end{document}